\documentclass[12pt, leqno, letterpaper]{amsart}

\usepackage{graphicx}    
\usepackage{bbm, xcolor}

\usepackage[margin=1in]{geometry}

\usepackage{amsmath,amsthm,amsfonts,amssymb}

\newtheorem{theorem}{Theorem}[section]

\newtheorem{lemma}{Lemma}[section]
\newtheorem{proposition}{Proposition}[section]

\numberwithin{figure}{section}

\definecolor{OliveGreen}{rgb}{0,0.6,0}

\begin{document}

\title[On the distribution of the time-integral of the gBM]
{On the distribution of the time-integral of the geometric Brownian motion}

\author{P\'eter N\'andori}
\address
{Department of Mathematical Sciences, Yeshiva University, New York, NY 10016}
\email{peter.nandori@yu.edu
}

\author{Dan Pirjol}
\address
{School of Business, Stevens Institute of Technology, Hoboken, NJ 07030}
\email{
dpirjol@stevens.com}

\date{August 2021}


\begin{abstract}
We study the numerical evaluation of several functions appearing in the
small time expansion of the distribution of the time-integral of the geometric Brownian motion 
as well as its joint distribution with the terminal value of the underlying Brownian motion.
A precise evaluation of these distributions is relevant for the simulation of stochastic volatility
models with log-normally distributed volatility, and Asian option pricing in the Black-Scholes model.
We derive series expansions for these distributions, which can be used for numerical evaluations.
Using tools from complex analysis, we 
determine the convergence radius and large order asymptotics of the coefficients in these expansions. We construct an efficient numerical approximation of the joint distribution of the time-integral of the gBM and its terminal value, and illustrate its application to Asian option pricing in the Black-Scholes model.
\end{abstract}

\maketitle

\baselineskip18pt

\section{Introduction}

The distributional properties of the time integral of the geometric Brownian motion (gBM)
appear in many problems of applied probability, actuarial science and mathematical finance. 
The probability distribution function is not known in closed form, 
although several integral representations have been given in the literature, as we describe next.

A closed form result for the joint distribution of the time integral 
of a gBM and its terminal value was obtained by Yor \cite{Yor}.
Denoting $A_t^{(\mu)} = \int_0^t e^{2(W_s+\mu s )} ds$, where $W_t$ is a standard Brownian
motion, this result reads
\begin{equation}\label{Yor0}
\mathbb{P}\left( \frac{1}{t} A_t^{(\mu)} \in da, W_t + \mu t \in dx \right) = 
e^{\mu x -\frac12 \mu^2 t}
e^{-\frac{1+e^{2x}}{2at}}
\theta_{\frac{e^x}{at}}(t) \frac{da dx}{a},
\end{equation}
where $\theta_r(t)$ is the Hartman-Watson integral defined by
\begin{equation}\label{HWthetadef}
\theta_r(t) = \frac{r}{\sqrt{2\pi^3 t}} e^{\frac{\pi^2}{2t}}
\int_0^\infty e^{-\frac{\xi^2}{2t}} e^{-r\cosh \xi} \sinh \xi \sin \frac{\pi\xi}{t} d\xi\,.
\end{equation}
and $r$ is a positive parameter.
This is related to the Hartman-Watson distribution \cite{Hartman} which is defined 
by its probability density function
$f_r(t) = \frac{\theta_r(t)}{\int \theta_r(s) ds}$.

Integrating (\ref{Yor0}) over $x$ gives an expression for the distribution of the time-average of the gBM as a double integral. For certain values of $\mu$ this can be reduced to a single integral
\cite{Dufresne2001}. Chapter 4 in \cite{Matsumoto2005} gives a detailed overview of the cases for which such simplifications are available.

The numerical evaluation of the integral (\ref{HWthetadef}) is challenging for small values of $t$
\cite{Barrieu}, which motivated the search for alternative numerical and analytical approximations. 
We summarize below limit theorems in the literature for the joint distribution of the time integral
of the gBM and its terminal value.

An asymptotic expansion for the Hartman-Watson integral
$\theta_{\rho/t}(t)$ as $t\to 0$ at fixed $\rho$ was derived in \cite{HWpaper}. 
The leading order of this expansion is
\begin{equation}\label{HWexp0}
\theta_{\rho/t} ( t ) = \frac{1}{2\pi t} 
e^{-\frac{1}{t}(F(\rho)-\frac{\pi^2}{2})}
G(\rho) + O(t^0)
\end{equation}
and depends on two real functions 
$F(\rho),G(\rho)$. 
The exact definition of $F$ and $G$ is technical and so we postpone it to Section 
\ref{sec:4.2}, see formulas  (\ref{Fsol}) and (\ref{Gsol}).

{Upon substitution into (\ref{Yor0}), the expansion 
(\ref{HWexp0}) gives the  leading $t\to 0$ 
asymptotics of the joint distribution of the time-integral of a gBM and its terminal value
\begin{equation}\label{Yorlimit}
\mathbb{P}\left( \frac{1}{t} A_t^{(\mu)} \in da, e^{W_t + \mu t} \in dv \right) = 
\frac{1}{2\pi t}  v^{\mu} e^{-\frac12\mu^2 t}
G(v/a) e^{-\frac{1}{t} I(a,v)} \frac{da dv}{a v} +  O(t^0)
\end{equation}
with
\begin{equation}
\label{eq:FI}
I(a,v) = \frac{1+v^2}{2a} + F\left( \frac{v}{a} \right) - \frac{\pi^2}{2} \,.
\end{equation}
The function $I(a,v)$ appears also in the asymptotic expansion of the log-normal SABR stochastic volatility model discretized in time, 
as the number of time steps $n\to \infty$
under appropriate rescaling of the model parameters, see Sec.~3.2 in \cite{PES}.

The application of Laplace asymptotic methods gives an expansion for the distribution function of the
time integral of the gBM, see Proposition 4 in \cite{HWpaper}
\begin{equation}\label{pdfA}
\mathbb{P}\left( \frac{1}{t} A_t^{(\mu)} \in da \right)  = \frac{1}{\sqrt{2\pi t}}
g(a,\mu) e^{-\frac{1}{t} J(a)} \frac{da}{a} ( 1 + O(t))
\end{equation}
where $g(a,\mu)$ is given in Proposition 4 of \cite{HWpaper}
and $J(a) = \inf_{v\geq 0} I(a,v)$. 

\subsection{Literature review on the small time asymptotics of the distribution of $(A_t^{(\mu)},W_t)$} 

In order to put these results into context, we summarize briefly the known short maturity 
asymptotics of the continuous time average of the geometric Brownian motion $a_T := \frac{1}{T}
\int_0^T e^{2(W_s+\mu s)} ds$. Dufresne \cite{Dufresne2004} proved convergence in law
of this quantity to a normal distribution, see Theorem 2.1 in \cite{Dufresne2004}
\begin{equation}\label{CLTa}
\lim_{T\to 0} \frac{a_T - 1}{\sqrt{T}} \stackrel{d}{=} N\left( 0, \frac43 \right)\,.
\end{equation}
This can be formulated equivalently as a log-normal limiting distribution in the same $T\to 0$ limit, see Theorem 2.2 in \cite{Dufresne2004}. The log-normal
approximation for the time-integral of the gBM is widely used in financial practice, see
\cite{Dufresne2004} for a literature review. When applied to the pricing of Asian options, this
is known as the L\'evy approximation \cite{Levy}.

The log-normal approximation is also used for the conditional
distribution of $A_t^{(\mu)}$ at given terminal value of the gBM, as discussed for example in \cite{McGhee,McGhee2011}. This conditional distribution appears in conditional Monte Carlo 
simulations of stochastic volatility models \cite{Willard1997}, see \cite{Brignone} for an overview
and recent results.

The result (\ref{CLTa}) can be immediately extended to the joint distribution of $(a_T,z_T)$ with
$z_T=e^{W_T+\mu T}$. We state it without proof which can be given along the lines of the 
proof of Theorem 2.1 in \cite{Dufresne2004}.

 \begin{proposition}\label{prop:CLT}
 Denote $\hat x_T = ( a_T, z_T)$ with $z_T := e^{W_T+\mu T}$.
 
As $T\to 0$ we have:
\begin{equation}\label{CLT2}
\lim_{T\to 0} \frac{\hat x_T- \hat 1}{\sqrt{T}} \stackrel{d}{=} N(0,\hat \Sigma)
\end{equation}
where $\hat 1 = (1,1)$ and 
\begin{equation}
\hat \Sigma = \left( 
\begin{array}{cc}
4/3 & 1 \\
1 & 1 \\
\end{array}
\right)
\end{equation}
\end{proposition}

The covariance matrix $\hat \Sigma$ is computed as 
(with $W_{s T} = \sqrt{T} B_s$ and $0\leq s \leq 1$)
\begin{eqnarray}
&& \Sigma_{11} = \mathbb{E}\Big[\Big(\int_0^1 2 B_s ds \Big)^2\Big] = 
4 \int_0^1 \int_0^1 \mbox{min}(s,u) ds du = \frac43 \\ 
&& \Sigma_{12} = \mathbb{E}\Big[\Big(\int_0^1 2 B_s ds \Big)(B_1)\Big] 
= 2\int_0^1 s ds = 1 \\ 
&& \Sigma_{22} = \mathbb{E}[(B_1)^2] = 1 \,.
\end{eqnarray}

In other words, as $T\to 0$, the joint distribution of $(a_T,z_T)$ approaches a bivariate normal
distribution with correlation $\rho_{az} = \frac{\sqrt{3}}{2}$. In analogy with
Theorem 2.2 in \cite{Dufresne2004}, this result can be formulated
equivalently as a bivariate log-normal limiting distribution.

Another asymptotic result in the $T\to 0$ limit is the large deviations property of the time-average
$a_T$. Using large deviations theory methods it was
showed in \cite{SIFIN} that, for a wide class of diffusions which includes the geometric Brownian 
motion, the time average of the diffusion satisfies a large deviations property in the
small time limit, see Theorem 2 in \cite{SIFIN}. For the particular case of a gBM this result reads, in our notations 
\begin{equation}\label{LDPa}
\mathbb{P}\left(\frac{1}{T} \int_0^T e^{2(W_s + \mu s)} ds \in da \right) = e^{-\frac{1}{4T} J_{\rm BS}(a )+o(1/T)}\,, \mbox{     as } T\to 0
\end{equation}
with $J_{\rm BS}(a)$ a rate function given in closed form in Proposition 12 of \cite{SIFIN}. For
convenience the result is reproduced below in Eq.~(\ref{JBS}). In contrast to the fluctuations result of
(\ref{CLTa}) which holds in a region of width $a-1 \sim \sqrt{T}$ around the central value $a=1$,
the large deviations result (\ref{LDPa}) holds for $a - 1 \sim O(1)$.

\subsection{Complete leading small time asymptotics}

The explicit asymptotic result  (\ref{Yorlimit}) for the small time asymptotics of the joint distribution
of $(A_t^{(\mu)},W_t)$ satisfies all these limiting results and improves on them by including 
all contributions up to terms of $O(t)$. It is instructive to see explicitly the connection to the
asymptotic results presented above. 

The exponential factor in (\ref{pdfA}) reproduces the large deviations result (\ref{LDPa}) by identifying 
$J(a) = \frac14 J_{\rm BS}(a)$. 
The emergence of the log-normal limit implied by Proposition ~\ref{prop:CLT} can be also seen explicitly by expanding the exponent in (\ref{Yorlimit}) around $a=v=1$ using the results of Propositions 2 and 3 in
\cite{HWpaper}.
Keeping the leading order terms in the expansion of $F(\rho)=\frac{\pi^2}{2}-1-\log \rho + \log^2\rho+O(\log^3\rho)$ in powers of $\log \rho$ gives the quadratic approximation
\begin{equation}\label{14}
I(a,v) = \frac32 \log^2 a - 3 \log a \log v + 2 \log^2 v + 
O(|\log^3 a| + |\log^3 v|) \,.
\end{equation}

The joint distribution (\ref{Yorlimit}) 
is thus approximated in a neighborhood of $(v,a)=(1,1)$ as
\begin{eqnarray}\label{joint}
\mathbb{P}\left( \frac{1}{t} A_t^{(\mu)} \in da, e^{W_t + \mu t} \in dv \right) = 
\frac{\sqrt3}{2\pi t} e^{-\frac{1}{2t}
[3\log^2 a - 6 \log a \log v + 4 \log^2 v +  \cdots ] }
C_{\rm LN}(a,v)
\frac{da dv}{a v} (1 + O(t)) \,,  \nonumber\\
\end{eqnarray}
where the $(\cdots )$ in the exponent are the error terms in (\ref{14}) and we defined
\begin{equation}
C_{\rm LN}(a,v) := v^\mu e^{-\frac12\mu^2 t} \frac{1}{\sqrt3} G(v/a)
= v^\mu e^{-\frac12\mu^2 t} ( 1 + O(\log (v/a))) \,.
\end{equation}

Up to the factor $C_{\rm LN}(a,v)$ and the error terms in the exponent, the expression (\ref{joint}) coincides with a log-normal
bivariate distribution in variables $(a,v)$. Recall the bivariate normal distribution of two random variables $N(X,Y;\sigma_X,\sigma_Y,\rho)$ with volatility and correlation parameters $\sigma_X,\sigma_Y,\rho$
\begin{equation}
\phi_2(x,y|\sigma_X,\sigma_Y,\rho) = \frac{1}{2\pi \sigma_X \sigma_Y \sqrt{1-\rho^2} }
e^{-\frac{1}{2(1-\rho^2)} \left( \frac{x^2}{\sigma_X^2} -2\rho \frac{x y}{\sigma_X\sigma_Y}
+ \frac{y^2}{\sigma_Y^2} \right) }\,.
\end{equation}
The density (\ref{joint}) corresponds to parameters
\begin{equation}
\sigma_X = \frac{2}{\sqrt3} \sqrt{t} \,, \sigma_Y = \sqrt{t}\,, \rho = \frac{\sqrt3}{2}\,.
\end{equation}
These parameters correspond precisely to the covariance matrix $\hat \Sigma$ 
appearing in the fluctuations result of Proposition \ref{prop:CLT}.

Improving on the log-normal approximation (\ref{joint}) requires adding higher order terms in the coefficient $C_{\rm LN}(a,v)$ and in the exponent. 
All these corrections are expressed in terms of $F(\rho),G(\rho)$ which can be evaluated by expansion around $\rho=1$. 
In this paper we present series expansions of $F(\rho),G(\rho)$ in powers of 
$\log\rho$ that can be used for such evaluation, and we study in detail their properties. 
We determine their convergence radii and coefficient asymptotics in Proposition \ref{prop:FG}.
We note that alternative series expansions which improve on the log-normal approximation have
been considered in the literature (in the context of the expansion of the density of $A_t^{(\mu)}$), based on Gram-Charlier series \cite{GramCharlier}.

Although closed form results are available for the functions $F(\rho),G(\rho)$, they are inefficient for practical use because of the need to solve one non-linear equation for each evaluation. The series expansions for $F(\rho),G(\rho)$ allow much faster and reliable numerical evaluation, inside their convergence domain.
In Section \ref{sec:5} we construct efficient numerical approximations for $F(\rho),G(\rho)$, based on the series expansions derived here within the convergence domain, which are easy to use in numerical applications. This gives a numerical approximation for the joint distribution of the time integral of the gBM and its terminal value, which appears in many problems of mathematical finance, such as the simulation of the SABR model \cite{Cai}
and option pricing in the Hull-White model \cite{HWmodel}.
As an application we discuss the numerical pricing of Asian options in the Black-Scholes model, and demonstrate good agreement with standard benchmark cases used in the literature.


\section{Preliminaries}

The mathematics of the series expansions for the functions $J(a),F(\rho),G(\rho)$ 
involves the study of the analyticity domain and
of the singularities of a function defined through the inverse of a complex entire function.
Such functions appear in several problems of applied probability, for example in 
applications of large deviations theory to the short maturity asymptotics of
option prices in stochastic volatility models \cite{Forde2009} and Asian 
options in the local volatility model \cite{SIFIN,AsianFwd}.
A similar problem is encountered in the evaluation of functions
defined by the Legendre transform of a function
\begin{equation}
F(x) = \inf_\theta (\theta x - f(\theta) ) = \theta_* x - f(\theta_*)
\end{equation}
where $\theta_*$ is the solution of the equation $f'(\theta) = x$.
If the optimizer $\theta_*$ can be found exactly at one particular point 
$x_0$, then expansion of $f(\theta)$ around $\theta_*$ 
gives a series expansion
for $F(x)$ in $x-x_0$ which can be used for numerical evaluation
of the Legendre transform.

For definiteness, consider the evaluation of the rate function $J_{\rm BS}$ 
appearing in the short maturity asymptotics of the distribution of the time-average of the gBM. This determines also the leading short maturity asymptotics of Asian options in the Black-Scholes model\footnote{More precisely, in the short maturity limit $T\to 0$ the prices of out-of-the-money Asian options 
approach zero at an exponential rate $C(K,T) \sim e^{-\frac{1}{T} J_{\rm BS}(K/S_0)}$,
determined by the rate function $J_{\rm BS}(x)$.}. 
The rate function can be computed
explicitly \cite{SIFIN} and is given by $J_{\rm BS}(x): (0,\infty) \to \mathbb{R}_+$, defined by
\begin{eqnarray}\label{JBS}
J_{\rm BS}(x)= \left\{
\begin{array}{cc}
\frac12\xi^2 - \xi \tanh \left( \frac{\xi}{2}\right) \,, & x\geq 1 \\
\zeta \tan \left(\frac{\zeta}{2}\right) - \frac12 \zeta^2 \,, & x \leq 1, \\
\end{array}
\right.
\end{eqnarray}
where $\xi$ is the solution of the equation
\begin{equation}\label{xidef}
\frac{\sinh \xi}{\xi} = x 
\end{equation}
and $\zeta$ is the solution in $[0, \pi]$ of the equation
\begin{equation}
\frac{\sin \zeta}{\zeta} = x \,.
\end{equation} 
See Proposition 12 in \cite{SIFIN}. 
The function $J_{\rm BS}(x)$ vanishes at 
$x=1$ and can be expanded in a series of powers of $x-1$ and $\log x$, see 
\cite{SIFIN}\footnote{Similar expansions were given for the rate functions of forward
start Asian options in \cite{AsianFwd}.}.  In practical applications one is interested in the
range and rate of convergence of these series expansions. 
We answer these questions in Proposition~\ref{prop:J}.

Section~\ref{sec:2} studies the analyticity properties of the complex inverse of the 
function determining $\xi, \zeta$ in (\ref{JBS}). These properties are used in Section~\ref{sec:3} to
study the singularity structure of the rate function $J_{BS}(x)$ and of the functions
$F(\rho),G(\rho)$ appearing in Eq.~(\ref{HWexp0}).
The nature of the singularity on the circle of
convergence determines the large-order asymptotics of the coefficients
of the series expansions of a complex function around a regular point. This is made precise by the transfer results (see Flajolet and Odlyzko \cite{Flajolet1988,Flajolet2008}) which relate the expansions around the singular point to the large-order asymptotics of the expansion coefficients. We apply these results in Section \ref{sec:3} to derive the leading large-order asymptotics of the expansion coefficients for the functions $J_{\rm BS}(x)$ and $F(\rho),G(\rho)$ around $x=\rho=1$,
see Proposition~\ref{prop:FG}.
Besides the rigorous proofs, numerical tests also confirm a general decreasing
trend of the error terms of the leading order asymptotics for the expansion coefficients.

\section{On the complex inverse of $x(\xi) = \frac{\sinh \xi}{\xi}$}
\label{sec:2}

An important role in the study of the expansions considered in this paper
is played by the singularities of a function defined through the inverse of another function, as in the
rate function $J_{\rm BS}(x)$, see Eq.~(\ref{JBS}). The main issues appear already in the 
study of the functions $\xi(x),\zeta(x)$. For this reason, we discuss the analytic structure of these 
functions in details in this section.

Let 
\begin{equation}\label{1}
f(\xi) = \frac{\sinh \xi}{\xi } \,.
\end{equation}
Clearly, $f$ is an entire function with Taylor series around zero given by 
\begin{equation}
\label{fexp}
f(\xi) = \sum_{n=0}^{\infty} \frac{\xi^{2n}}{(2n+1)!}.
\end{equation}
We want to make sense of the inverse of $f$, restricted to some complex domains.
First observe that the
equation $\sinh \xi= \xi$ has a unique real solution $\xi=0$
and infinitely many complex solutions. 
This shows that $f$ is not
invertible on the entire complex plane. In this section
we discuss how to make sense of an inverse. Although the results
will not surprise experts in complex analysis,
this question has not been discussed in the present context before.

\subsection{Local analysis}

Let us now restrict $f$ to a small
neighborhood of the origin. Then 
the range of $f$ is in a small neighborhood of $1$. 
Clearly, $f$ is not locally
invertible as $f$ is an even function, that is $f(\xi) = f(-\xi)$, 
and thus is not injective. 
A natural way to make the function at least locally invertible is to define 
\begin{equation}
\label{gexp}
g(z) = \sum_{n=0}^{\infty} \frac{z^n}{(2n+1)!} \,.
\end{equation}
In light of \eqref{fexp}, we have 
$g(z) = f(\sqrt{z}) = \sinh (\sqrt z) / \sqrt z$. 
Here, and in the sequel, unless stated otherwise, the notation 
$\sqrt z$ will only be used in the argument of even functions so as both values of 
$\sqrt z$ give the same result.
Recall the following simple fact:
if $g'(z) \neq 0$, then there is a small neighborhood $D$ of $z$
so that $g: D \to g(D)$ is biholomorphic (see for example
Theorem 3.4.1 in \cite{Simon}).
Since $g'(0) \neq 0$, $g$ is locally invertible at zero. 
We also note that $g$ is an entire function as its Taylor series
converges everywhere. We will show next that $g: \mathbb C \mapsto \mathbb C$ is surjective. 
One says that an entire function $f$ has finite order $\rho$ if constants $A,B$
can be found such that $|f(z)| \leq A e^{B |z|^\rho}$ for all $z\in \mathbb{C}$. 
See for example Chapter 5 in \cite{Stein2003} and Section 1.3 in \cite{Levin}
for a detailed discussion of entire functions of finite order.

The order of the entire function 
$g(z)=\sum_{n=0}^\infty a_n z^n$ can be expressed in terms of the coefficients 
of its Taylor expansion as, see Theorem 2 in Sec.~1.3 of \cite{Levin}
$$
\rho = \limsup \frac{n \log n}{- \log |a_n|} \,.
$$
Using the expansion (\ref{gexp}) it follows that the function $g$ has order $\rho=\frac12$. 
Any entire function with order $< 1$ is surjective, see Section 5.1 of \cite{Levin}.
Since in our case $\rho = 1/2$, $g: \mathbb C \mapsto \mathbb C$ is surjective.

We study next the critical points of $f$ and $g$ which are given by the zeros of $f'$ and $g'$, respectively. To identify the zeros
of $g'$, let us write
$0=\eta_0 < \eta_1 < \eta_2< ...$ for the 
non-negative solutions of the equation $\tan \eta = \eta$ (and note that all real 
solutions are given by $\eta_k, k \in \mathbb Z$ with $\eta_{-k} = - \eta_k$).
We also write 
\begin{equation}
\label{def:const}
\omega_k = \frac{\sin (\eta_k)}{\eta_k}, \quad \xi_k = i \eta_k, \quad
z_k = -\eta_k^2
\end{equation} 
(see Table \ref{Tab:1} for the numerical values for $k=1,...,5$).
Now we have

\begin{lemma}
\label{lem1}
The zero sets of $f'$ and $g'$ are given by
$$
\{ \xi: f'(\xi) = 0\} = \{ \xi_k, k \in \mathbb Z\}
$$
$$
\{ z: g'(z) = 0\} = \{ z_k, k \in \mathbb Z_+\}
$$
\end{lemma}

\begin{proof}
If $\xi \neq 0$, then the equation $f'(\xi) = 0$ reduces to $\xi = \tanh (\xi)$,
whose solutions are exactly the numbers $\xi_k$, $k \neq 0$. In case $\xi = 0$,
we see by \eqref{fexp} that $f'(0) = 0$. The first statement follows.
For the second statement, observe that $g'(z) = f'(\sqrt z) /(2 \sqrt{z})$
and $g'(0)  = 1/6 \neq 0$. 
\end{proof}

\begin{lemma}
\label{lem2}
For every $k \in \mathbb Z_+$, $g''(z_k) \neq 0$.
\end{lemma}

\begin{proof}
Since $g$ maps $\mathbb R$ to $\mathbb R$, in this proof we will denote by 
$\sqrt z$ for the square root
function from $\mathbb R_+$ to $\mathbb R_+$.
By Lemma \ref{lem1}, $g'(z_k) = 0$, that is $\tanh (\sqrt z_k) = \sqrt z_k$. Next, 
we compute
$$
g''(z) = \frac{(z+3)\sinh (\sqrt z) - 3 \sqrt z \cosh (\sqrt z)}
{4 z^2\sqrt z }
$$
Assuming $g''(z_k) = 0$, we obtain 
$\tanh (\sqrt z_k) = \frac{3 \sqrt z_k}{z_k+3}$ which
is a contradiction with 
$\tanh (\sqrt z_k) = \sqrt z_k$ as $z_k \neq 0$.
\end{proof}

We note that $\omega_k \in \mathbb R$ for all 
$k \geq 1$. Furthermore, 
$\omega_1 \approx -0.217234$, and $|\omega_k| < |\omega_1|$ for all $k\geq 2$.

Let $B_{\rho}(A)$ be the open $\rho$-neighborhood of the set $A \subset \mathbb C$
and $B_{\rho}(w) = B_{\rho}(\{ \omega \})$. Now we are ready to state the analyticity of the local inverse of $g$
on the region that is most important for our applications.
\begin{theorem}
\label{thm1}
There is a function $h$, analytic on $B_{1-\omega_1}(1)$ so that $h(1) = 0$ and 
$g(h(\omega)) = \omega$ for all $w \in B_{1-\omega_1}(1)$.
\end{theorem}

\begin{proof}
Fix an arbitrary $\varepsilon >0$. It suffices to show that $h$ can be defined analytically
on ${\bf B}:=B_{1-\omega_1 - \varepsilon}(1)$.
First observe that the real interval $I = [z_1 + \varepsilon_-, \varepsilon_+]$ is being mapped,
under the function $g$, to the real interval $g(I) = [\omega_1 + \varepsilon, 1+\varepsilon]$
for suitable $\varepsilon_{\pm}$. Furthermore, the derivative of $g$ is non-zero
on $I$. Thus restricting the domain of $g$ to a small neighborhood of $I$, 
the restricted function provides a bijection between a neighborhood of
$I$ and a neighborhood of $g(I)$. More precisely, since the derivative
of $g$ on $I$ is non-zero,
for every $\omega \in g(I)$ there is some $\delta_{\omega}$ so that
$g$ can be analytically inverted on $B_{\delta_\omega}(\omega)$. By compactness of $g(I)$,
we can choose a finite cover of $g(I)$ by the sets $B_{\delta_\omega}(\omega)$ showing that
the inverse function $h$ is well defined on a small neighborhood 
$B_{\delta'}(g(I))$ of $g(I)$.

Next we extend the function $h$ analytically to the ball $\bf B$.
This can be done by a standard procedure. 
Recall that a point $a \in \mathbb C$ is an asymptotic value of $g$ if there is a curve
$\Gamma$ tending to $\infty$ so that $g(z) \rightarrow a$ az $z \rightarrow \infty$
along $\Gamma$.
Furthermore, if $g'(c) = 0$, then $g(c)$ is a critical value. The singularity set
is defined as the union of asymptotic values and critical values. It is known
that for any open set $D$ that is disjoint to the singularity set, the map
$g:g^{-1}(D) \rightarrow D$ is a covering (see \cite{N}).

In our case, we have already identified the set of critical values $\{ \omega_k, k \geq 1\}$.
An elementary computation shows that the only asymptotic value is $\omega = 0$.
Now let $D = {\bf B } \setminus \overline{ B_{\delta'/2}({g(I)})}$ (where $\bar B$ is the closure of $B$). Then $D$ is open and so
$g:g^{-1}(D) \rightarrow D$ is a covering. This implies that along any simple polygonal
arc $\gamma$
starting from $B_{\delta'/2}(I)$ and staying inside $\bf B$, $h$ has an analytic extension. 
Indeed, by the covering map property, every $\omega \in \gamma$
has a neighborhood $U_\omega$ so that $g^{-1}(U_\omega)$
 is a union of disjoint open sets, each of which is mapped homeomorphically onto U by $g$. Let us now choose a finite cover of $\gamma$ by neighborhoods 
$U_{\omega_1}, ..., U_{\omega_n}$ so that $U_{\omega_1} \cap B_{\delta'/2}(I) \neq \emptyset$
and  $U_{\omega_{k-1}} \cap U_{\omega_{k}} \neq \emptyset$ and by induction on $k$
we can extend $h$ analytically to $\cup U_{\omega_k}$. Observing that $D$ is simply 
connected, the result now follows from the monodromy theorem.

\end{proof}

\subsection{Global analysis}
\label{sec:global}

Theorem \ref{thm1} tells us that the inverse 
$h(\omega)$ can be defined with $h(1) = 0$. Furthermore,
$h$ is analytic and the Taylor expansion 
\begin{equation}
\label{invexpansion}
h(\omega) = \sum_{n=1}^{\infty} c_n (\omega-1)^n
\end{equation}
has a radius of convergence of at least $1-\omega_1$. 
We note that the inverse $g^{-1}$ can be defined in the
sense of a sheaf (also known as global analytic function, or Riemann surface). 
In this general sense,
$g^{-1}$ has a second order algebraic branch point
at $\omega_1 = g(z_1)$ by Lemma \ref{lem2}. Thus the function
$h$ cannot be analytically extended to $\omega_1$ and so the 
power series \eqref{invexpansion} has radius of 
convergence exactly $1-\omega_1$. We will be mostly interested 
in this expansion (see Section \ref{sec:exp}).
However, for completeness we first discuss how to 
visualize the sheaf $g^{-1}$.

First, observe that $h$ can be analytically extended to 
$\mathbb C \setminus \Gamma_1$, where $\Gamma_1$ is any 
curve starting at $\omega_1$ and tending to $\infty$ (for example
$\Gamma_1 = \omega_1 + i\mathbb R_{\geq 0}$; however as we will see later, sometimes
other choices are more convenient). Indeed, 
replacing $D$ by 
$$\mathbb C \setminus (\overline{ B_{\delta'/2}({g(I)})} \cup
\{ x+iy: |x - \omega_1| \leq \varepsilon, y \geq - \varepsilon\}
)$$
in the last step of the proof of Theorem \ref{thm1} and noting
that this domain is simply connected, 
analyticity of $h$ follows (of course the expansion \eqref{invexpansion} is not valid on the larger domain).
The function $h$ can be identified with a chart of the Riemann
structure of the sheaf $g^{-1}$ (in this
identification, its domain is sometimes referred to as
Riemann sheet). 

Next, we can define another analytic 
function $h_2$ on $\mathbb C \setminus (\Gamma_1 \cup \Gamma_2)$,
where $\Gamma_2 = \omega_2 + i\mathbb R_{\geq 0}$.
Since the function $z \mapsto g(z) - \omega_1$
has a zero of order $2$ at $z_1$, we can extend $h$ analytically
from 
$$\{\omega_1 + re^{i \theta}: r < \varepsilon, \theta \in (\pi/2, 3\pi/4) \}$$
to 
$$\{\omega_1 + re^{i \theta}: r < \varepsilon, \theta \notin [3\pi /4, \pi] \}.$$
Let the resulting function be denoted by $\tilde h$. Clearly,
$\tilde h$ takes values in a small neighborhood of $z_1$
and since there is a branch point at $\omega_1$, $\tilde h (\omega) \neq h(\omega)$ for 
$\omega = \omega_1 + t$ with $t \in \mathbb R_+$, $t < \varepsilon$. 
Now observe that 
$g((z_2,{z_1}]) = [\omega_1, \omega_2)$. Since
the branch point is of order $2$, we 
conclude that $ \tilde h(\omega_1 + t)$ for $t$ as above 
has to be real, and slightly bigger than $z_1$. 
Thus we can extend $\tilde h$ analytically to 
$B_{\delta}(g(I_2))$ where $I_2 = [z_2 + \varepsilon_- , z_1 - \varepsilon_+]$
and $g(I_2) = [\omega_1 + \varepsilon/2, \omega_2 - \varepsilon/2]$. Now restrict $\tilde h$
to 
$$\{\omega_1 + re^{i \theta}: r < \varepsilon, \theta \notin [\pi /2, \pi] \} \cup B_{\delta}(g(I_2))$$
and denote this restriced function by $h_2$. Then $h_2$
can be extended analytically to 
$\mathbb C \setminus (\Gamma_1 \cup \Gamma_2)$ by the monodromy theorem
as in the proof of Theorem \ref{thm1}.

Furthermore, we can "glue together" the domains of $h_1$ and $h_2$
along $\Gamma$ since
$\lim_{x \rightarrow \omega_1 \pm} h_1(x + yi) = \lim_{x \rightarrow \omega_1 \mp} h_2(x + yi) $ and the gluing 
is analytic (see the analyticity of $\tilde h$ and the fact that the branch point is algebraic of order $2$). Since
$g$ preserves the real line and all branch points are of second order, 
this construction can be continued inductively, i.e. we can 
define the Riemann sheet $k \geq 2$
by $h_k$ and glue together with the sheets $k \pm 1$ along $\Gamma_{k\pm 1}$.
This gives a complete elementary description of the sheaf $g^{-1}$.

\subsection{Visualizing the sheaf}
For better visualization of the sheaf constructed in the 
previous section, we include some figures. We start with listing the numerical 
values of the first 
few critical points $\eta_k$, $\eta_k^2$, $\omega_k$ and $|i\pi + \log |\omega_k||$
(the latter will be needed later)
 in Table \ref{Tab:1}.

\begin{table}
\caption{\label{Tab:1}
Numerical values of the first few numbers of interest as defined in \eqref{def:const}.
}
\begin{center}
\begin{tabular}{|c|ccc|c|}
\hline 
$k$ & $\eta_k$ & $z_k=-\eta_k^2$ & $\omega_k = \sin(\eta_k) /\eta_k $ & $|i\pi + \log |\omega_k||$ \\
\hline \hline
1 & 4.4934 & -20.19 & -0.2172 & 3.4929 \\
2 & 7.7252 & -59.68 &  0.1284 & 2.0528  \\
3 & 10.9041 & -118.90  & -0.0913 &  3.9494 \\
4 & 14.0662 & -197.86  &  0.0709 &  2.6463 \\
5 & 17.2208 & -296.55  &  -0.0580 &  4.2402 \\
\hline
\end{tabular}
\end{center}
\end{table}

\begin{figure}
\centering
\includegraphics[width=3.0in]{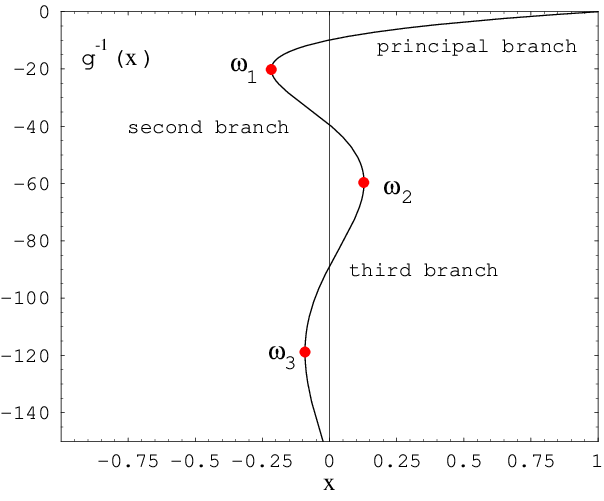}
\caption{Branches of $g^{-1}(x)$ for $x\in (\omega_1,1)$ along the real axis.}
\label{Fig:branches}
\end{figure}

Given $\omega$, there may be infinitely many $z$'s so that $g(z) = \omega$. 
A plot of some of these $z$'s
for $\omega \in (\omega_1, 1)$ are shown in 
Fig.~\ref{Fig:branches} (note that both $\omega$ and $z$ are real).
The first branch corresponds to
$z \in (z_1, \infty)$ when $\omega > \omega_1$.
On this branch $h_1(1)=0$. We call this Riemann sheet the principal branch. The second branch corresponds to 
$z \in (z_2, z_1)$, where
$\omega \in (\omega_2,\omega_1)$, 
and the $k$-th branch corresponds to $z \in (z_k, z_{k-1})$, where
$\omega$ is in between $\omega_{k-1}$ and $\omega_k$.

\begin{figure}
\centering
\includegraphics[width=3.6in]{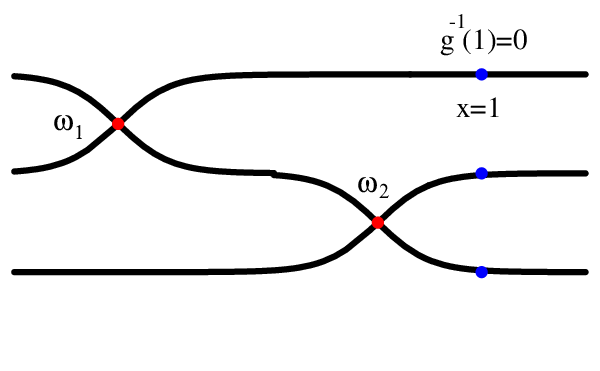}
\caption{
Schematic representation of the first three Riemann sheets
of the function $g^{-1}(x)$. 
The blue dots show $g^{-1}(1)$ on the various sheets. 
On the main sheet $g^{-1}(1)=0$.
This is joined to the second sheet at the branch cut $\Gamma_1$ starting
at the branch point $\omega_1\approx-0.217$. This is joined to the third
sheet at the branch cut $\Gamma_2$ starting at the branch point $\omega_2\approx 0.128$,
and so on. There is an infinite number of Riemann sheets joined successively
at the branch points $\omega_k\to 0$ which accumulate towards zero as $k\to \infty$.} 
\label{Fig:Riemann}
\end{figure}

Figure \ref{Fig:Riemann} is a schematic representation of the first three Riemann
sheets. This figure illustrates the apparently counterintuitive result that, although the branch point at $\omega_2=0.1284$
is closer to $\omega=1$ than the branch point at $\omega_1=-0.2172$, 
it is the latter which determines the radius of 
convergence of the series expansion around $\omega=1$, rather than the former. 
The explanation for this fact is that the branch point at $\omega_2$ is not ``visible'' on the
main Riemann sheet, and thus the branch point at $\omega_1$ is the dominant singularity
around $\omega=1$.

\subsection{Properties of the Taylor series}
\label{sec:exp}

Next, we study the
coefficients $c_n$ of the Taylor expansion \eqref{invexpansion}. The first few coefficients
are easily computed by the Lagrange inversion theorem using
\eqref{gexp} and Theorem \ref{thm1}. In particular, we obtain

\begin{equation}\label{1inv}
h(\omega) = 6(\omega-1) - \frac95 (\omega-1)^2 + \frac{144}{175} (\omega-1)^3 + 
O((\omega-1)^4)
\end{equation}

Writing $e^y= \omega$ we can also expand $h$ as a series in $y=\log \omega$
(here, $\log \omega$ means the principal value of the logarithm).
In applications to Asian options pricing \cite{SIFIN}, the parameter $\omega$ 
is related to the option strike. In practice, it is usual to expand implied volatilities
in log-strike, which corresponds to expanding in $y=\log\omega$. See for example \cite{Gatheral}.

Comparing the first few derivatives with \eqref{1inv}, the expansion (\ref{1inv})
becomes 
\begin{eqnarray}\label{1invlog}
h(e^y) &=& \sum_{n=1}^{\infty} d_n y^n \\
&=& 6y + \frac65 y^2 + \frac{4}{175} y^3 - \frac{2}{175} y^4 + O(y^5). \nonumber
\end{eqnarray}

Explicitly computing high order terms of these expansions is difficult
(cf. \cite{MF53}, pages 411-413). However, we can study
the asymptotic properties of the sequences $c_n$ and $d_n$.

\begin{proposition}\label{prop:2}
(i) The 
coefficients $c_n$ of the Taylor expansion \eqref{invexpansion} satisfy
\begin{equation}\label{cnasympt}
c_n = c_\infty \frac{1}{(1-\omega_1)^n} (-1)^n n^{-\frac32} + O(n^{-\frac52})\,,\quad
c_\infty =  - \eta_1 \sqrt{\frac{2(1-\omega_1)}{\pi |\ell''(\eta_1)|}} \approx -8.48671,
\end{equation}
where $\ell (\eta) = \frac{\sin\eta}{\eta}$.

(ii) The Taylor expansion (\ref{1invlog}) has a radius of convergence $\rho_x = |\log |\omega_1| + i\pi| \approx 3.49295$. 
Furthermore, the 
coefficients $d_n$ satisfy
\begin{equation}\label{dnasympt}
d_n = d_\infty \frac{1}{\rho_x^n} \cos\left[\theta_x \left(n-\frac12\right)\right]
n^{-\frac32} + O(n^{-\frac52})\,,\quad
d_\infty =  -2\eta_1 \sqrt{\frac{2(-\omega_1)\rho_x}{\pi |\ell''(\eta_1)|}}\approx -13.4011,
\end{equation}
where $\theta_x=\mbox{arg}(i\pi + \log |\omega_1|) \approx 2.02317$.

\end{proposition}

\begin{figure}
\centering
\includegraphics[width=2.6in]{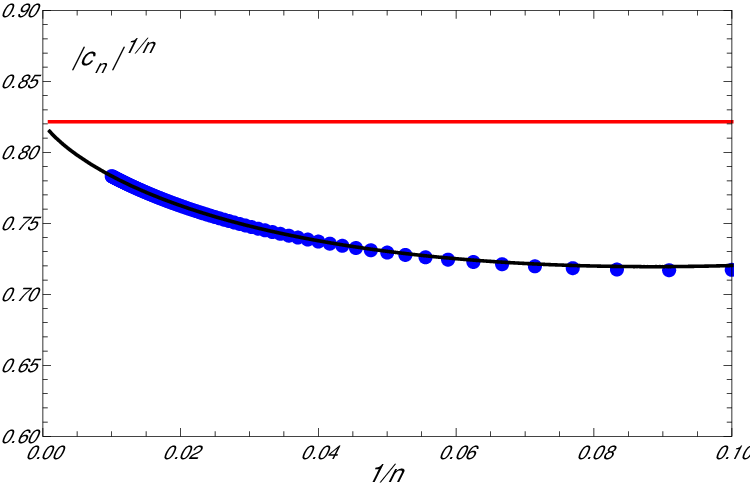}
\includegraphics[width=2.6in]{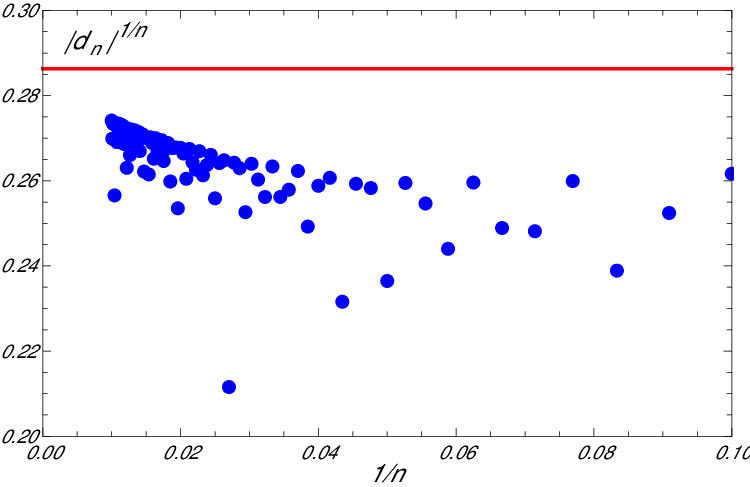}
\caption{Left: plot of $[c_n]^{1/n}$ vs $1/n$ for the
series (\ref{1inv}) using the first $100$ terms (blue dots).
The red line is at 
$1/(1-\omega_1)=0.821$. The black solid curve shows the leading 
$n\to \infty$ asymptotics for the coefficients $c_n$ (\ref{cnasympt}).
Right: plot of $[d_n]^{1/n}$ for the series (\ref{1invlog}) using the 
first $100$ terms (blue dots).
The red line is  at $1/\rho_x=0.286$ with $\rho_x=|\log |\omega_1| + i\pi|=3.49295$.}
\label{Fig:root}
\end{figure}

\begin{proof} (i) We first prove \eqref{cnasympt}.
The asymptotics of the coefficients $c_n$ for $n\to \infty$ depend on the
nature of the singularity on the circle of convergence. The leading asymptotics
of $h(\omega)$ around this singularity can be translated into asymptotics for
the coefficients of the series expansion using the "transfer results" (see Flajolet 
and Odlyzko \cite{Flajolet1988,Flajolet2008}). Theorem VI.1 and VI.3 of \cite{Flajolet2008}
imply that if $h$ is analytic on $\mathbb C \setminus \Gamma'$
where $\Gamma' = \mathbb R_{\leq \omega_1}$ and furthermore
\begin{equation}
\label{eq:exphomega_1}
h(\omega ) =  h(\omega_1) + C_1 (\omega - \omega_1)^{1/2} + C_2(\omega - \omega_1) + O((\omega- \omega_1)(\omega - \omega_1)^{1/2}) 
\end{equation}
as $\omega \rightarrow \omega_1$,
\footnote{Here $(\omega - \omega_1)^{1/2}$ means an analytic function on $B_{\varepsilon}(0) \setminus \mathbb R_{\leq 0}$
whose square is $(\omega - \omega_1)$}
then \eqref{cnasympt} holds with 
\begin{equation}
\label{abstractc}
c_{\infty} = - C_1 \frac{\sqrt{1-\omega_1}}{2 \sqrt \pi}.
\end{equation}
We defined here $h(\omega_1) = z_1$. 
With this definition,
$h$ is continuous at $\omega_1$ (of course there is no open neighborhood of $\omega_1$ where $h$ is defined).

Indeed, Theorem VI.1 of \cite{Flajolet2008} says that the Taylor coefficients $a_n$
of the function $(1-\omega)^{1/2} = \sum a_n \omega^n$ have the large $n$
asymptotics $a_n = \frac{n^{-3/2}}{\Gamma(-1/2)} + O(n^{-5/2})$. 
The term $C_2(\omega - \omega_1)$ is a polynomial and so it has a finite
Taylor expansion and by Theorem VI.3 of \cite{Flajolet2008}, the Taylor coefficients of $(1-\omega)(1-\omega)^{1/2}$ are
of $O(n^{-5/2})$. We obtain \eqref{abstractc} by rescaling the singularity from $1$ to $\omega_1$.

To finish the proof of \eqref{cnasympt}, it remains to verify \eqref{eq:exphomega_1}.
We will prove that $\omega_1$ is an algebraic branch point of order 2, 
meaning that near $\omega_1$, the function $h$ 
``behaves like a branch of  the square root function near $0$'' . 
This behavior results from the inversion of the function $f(\sqrt{z})$ around the
critical point $z_1$. 
It is known, see for example Theorem 3.5.1 in \cite{Simon}, that the inversion of a
non-constant analytic function $f(z)=w$ around a point $z_0$ with $f'(z_0)=0$
gives a multi-valued function which can be expanded in a convergent series in fractional powers of $w-w_0$ with $w_0=f(z_0)$, called a Puiseux series.

To make these statements mathematically precise, first note that $h$ can be analytically defined on 
$\mathbb C \setminus \Gamma'$
exactly as in Section \ref{sec:global} (replacing $\Gamma_1$ by $\Gamma'$. With a slight abuse of notation we also denote this function
by $h$.)
Then  by Theorem 3.5.1 in \cite{Simon}
there is some $\varepsilon >0$ and an 
analytic function $\phi$ on $D' := B_{\varepsilon}(0) \setminus \mathbb R_{\leq 0}$
so that $\phi^2(x) = x$ (i.e. $\phi$ is one branch of the square root function) and another
function $\psi$, analytic on $B_{\varepsilon}(0)$ with 
$\psi(0) = z_1$, $\psi'(0) \neq 0$
so that $h(\omega - \omega_1) = \psi(\phi(\omega))$
for all $\omega \in D'$.
Furthermore, the function $h$ has a Pusieux expansion on $D'$, that is

\begin{equation}
\label{hexpphi}
h(\omega) = z_1 + C_1 \phi(\omega - \omega_1) + C_2 (\omega - \omega_1) + O ((\omega - \omega_1)
\phi(\omega - \omega_1))
\end{equation}
where $C_1$ satisfies $C_1^2 = -8z_1/ f''(i \xi_1) $. Thus we obtain \eqref{eq:exphomega_1} which 
completes the proof of \eqref{cnasympt}.

(ii)
First we prove that the Taylor expansion (\ref{1invlog}) has a radius of convergence $ \rho_x := |\log |\omega_1| + i\pi|$.
First note that the right hand side of (\ref{1invlog}) converges to an analytic function $\tilde h$ 
if $|\Im y |< \pi$. Indeed, in this case,
$e^y \notin \mathbb R_-$ and so $\tilde h = h(e^y)$ is well defined as $h$ is analytic on $\mathbb C \setminus \Gamma'$.
Now we claim that $\tilde h$ can be analytically extended to $\tilde D = B_{\rho_x}(0)$. To prove this claim, let 
$\tilde D_+ = \tilde D \cap \{ y: \Im y > \pi \}$ and $\tilde D_- = \tilde D \cap \{ y: \Im y < - \pi\}$. 
Let us denote by $\mathbb H$ the upper half plane, that is $\mathbb H = \{ y \in \mathbb C: \Im y >0\}$.
Now observe that the image
of $\tilde D \cap \mathbb H$ under the exponential function
is disjoint to $\Gamma_+ : = \omega_1 + e^{i (2 \pi - \epsilon)} \mathbb R_{\geq 0}$ for $\epsilon$ small enough.
Thus defining $h$ with the branch cut $\Gamma_+$, $h(e^y)$ is well defined and analytic for $\tilde D \cap \mathbb H$. Likewise, we can analytically
extend $\tilde h$ to $\tilde D_-$ using the branch cut $\Gamma_- : = \omega_1 + e^{i  \epsilon} \mathbb R_{\geq 0}$ in the definition of $h$.
Clearly, there are singularities at $y_\pm = \pm i\pi
+ \log |\omega_1|$ as $e^{y_\pm} = \omega_1$. Thus we have verified that the right hand side of (\ref{1invlog})
converges to an analytic function $\tilde h$ on $B_{\rho_x}(0)$.

Next, we prove \eqref{dnasympt}.
The new feature here is the presence of two singularities on the
circle of convergence in the $y$ complex plane, with radius $\rho_x$. This situation is discussed in Chapter VI.5 of \cite{Flajolet2008}.
One feature specific to this situation is 
the oscillatory
pattern of the leading
asymptotics for the coefficients of the series expansion 
\eqref{1invlog}.

From \eqref{eq:exphomega_1}, we have for $y\to y_\pm$
\begin{eqnarray}
\tilde h(y) &=& z_1 + C_1 (e^y - e^{y_{\pm}})^{1/2}
+ C_2 (e^y - e^{y_{\pm}}) + O( (e^y - e^{y_{\pm}}) (e^y - e^{y_{\pm}})^{1/2})
\\
&=& z_1 + C_1 (y - {y_{\pm}})^{1/2} \sqrt{\omega_1}
+ C_2 (e^y - e^{y_{\pm}}) + O( (e^y - e^{y_{\pm}})  (y - {y_{\pm}})^{1/2}).
\nonumber
\end{eqnarray}
Applying again the transfer result, we obtain
\begin{equation}
d_n = C_1
[(y_+)^{\frac12-n} + (y_-)^{\frac12-n}] \frac{-1}{2\sqrt{\pi n^3}} + O(n^{-\frac52})
\end{equation}
which reproduces (\ref{dnasympt}) after substituting here $y_\pm = \pm i\pi + \log |\omega_1|$.

\end{proof}

We present next some numerical tests of these asymptotic results.
In Fig.~\ref{Fig:root} (left) we plot $|c_k|^{1/k}$ vs $1/k$ for the first 100 terms
of the series (\ref{1inv}) (left panel) and $|d_k|^{1/k}$ vs $1/k$ for the first 
100 terms of the series (\ref{1invlog}) (right panel).
The red lines in these plots are at $1/|1-\omega_1|$ (left) and $1/|\log |\omega_1|+i\pi|$ 
(right). The general trend as $k\to \infty$ is towards the limiting result shown by 
the red lines is consistent with Proposition \ref{prop:2}.

\section{Applications: Series expansions}
\label{sec:3}

We apply in this section the results of the previous section to two specific
expansions appearing in mathematical finance which were discussed in the Introduction.

\subsection{The rate function $J_{BS}(x)$}

The rate function $J_{\rm BS}$ appears in the small-time asymptotics for the 
density of the time average of a gBM, and determines the leading short maturity
asymptotics for Asian options in the Black-Scholes model. 
The explicit form of this function for real positive argument $x>0$
was given above in equation \eqref{JBS}.
Complexifying the function $J_{BS}(x)$ and denoting the complex
variable by $\omega$, $J_{\rm BS}$ can be written as
\begin{equation}
J_{\rm BS}(\omega) = \mathcal{J}(h(\omega)),
\end{equation}
where $h(\omega)$ is the function introduced in Theorem~\ref{thm1}, and we define 
$$\mathcal{J}(z)=\frac12 z - \sqrt{z} \tanh(\sqrt{z}/2).$$
Note that $\mathcal J$ is well defined as the function $z \mapsto {z} \tanh({z}/2)$ is even.

\begin{proposition}\label{prop:J}
(i)The series expansion of the rate function $J_{BS}(\omega)=\sum_{n=0}^\infty c_{J,n}(\omega-1)^n$ 
has convergence radius $R_J=1$. The first few terms are given by
\begin{equation}\label{JBSexpomega}
J_{\rm BS}(\omega) = \frac32 (\omega-1)^2 - \frac95 (\omega-1)^3 + \frac{333}{175} (\omega-1)^4
+ O((\omega-1)^5)
\end{equation}
and the $n\to \infty$ asymptotics of the coefficients $c_{J,n}$ are given by
\begin{equation}
\label{JBScoeffexp}
c_{J,n}= (-1)^n 2 +O((1-\omega_1)^{-n})\,.
\end{equation}

(ii) The series expansion $J_{BS}(e^y)=\sum_{n=0}^\infty d_{J,n} y^n$
around $y=0$ converges for $|y | < \rho_x = |i\pi + \log |\omega_1||\approx 3.49295$.
The first few terms are
\begin{equation}\label{JBSexplog}
J_{\rm BS}(e^y) = \frac32 y^2 - \frac{3}{10} y^3 + 
\frac{109}{1400} y^4 + O(y^5)\,.
\end{equation}
The large $n$ asymptotics of the coefficients is
\begin{equation}
\label{JBScoeffexp2}
d_{J,n} = d_J \rho_x^{-n} \cos\left[\theta_x\left(\frac32 - n\right)\right]
n^{-\frac52} + O(n^{-\frac72})\,.
\end{equation}
with $d_J \approx -23.4048$ and $\theta_x \approx 2.02317$.
\end{proposition}

\begin{proof}
(i) First,
observe that the function $\mathcal{J}$ has isolated singularities (simple poles)
along the real negative axis at $z_n = -(2n + 1)^2 \pi^2$ for $n \in \mathbb Z$, and is analytic everywhere else. 

The residues of $J(z)$ at the simple poles $z_n$ are 
$\mbox{Res } \mathcal{J}(z=z_n) = -4z_n$.
Recall that by Theorem 2.1, $h$ is analytic on $B_{1-\omega_1}(1)$. 
Thus for any $\omega \in B_{1-\omega_1}(1)$
with $h(\omega) \neq - (2n + 1)^2 \pi^2$, the 
function $J_{\rm BS}(\omega) $ is analytic at $\omega$. 
Noting that the only solution of the equation $h(\omega) = -(2n + 1)^2 \pi^2$ is $\omega = 0$, we conclude that
$J_{\rm BS}(x) $ is analytic on 
$B_{1-\omega_1}(1) \setminus \{ 0 \}$ and has a singularity at $0$.  
Thus $R_J =1$.

Next, we prove \eqref{JBScoeffexp}. In order to apply the transfer results (cf. Flajolet and Odlyzko)
we need to study the behavior of $J_{\rm BS}(\omega)$ around the singularity
on the circle of convergence. This is a simple pole at $\omega=0$, so it is sufficient
to compute its residue.

Recall that $\mbox{Res } \mathcal{J}(z=- \pi^2) = 4\pi^2$.
Thus, as $\omega \to 0$ we have
\begin{eqnarray}\label{JBSexp}
J_{BS}(\omega) = \mathcal{J}(h(\omega)) = 
\frac{4\pi^2}{h(\omega) + \pi^2} + \mbox{ a bounded term as } \omega \to 0 \,.
\end{eqnarray}

The Taylor expansion of $h(\omega)$ around $\omega=0$ is 
\begin{equation}
h(\omega)= h(0) + h'(0) \omega + O(\omega^2)
\end{equation}
where $h(0)=-\pi^2$ and 
\begin{equation}
h'(0) = \frac{dh}{d\omega}|_{\omega=0}=
\frac{1}{\frac{d}{dh} \frac{\sinh \sqrt{h}}{\sqrt{h}}}|_{h=-\pi^2} = 2\pi^2\,.
\end{equation}

Substituting into (\ref{JBSexp}) gives
\begin{equation}
J_{BS}(\omega) = \frac{2}{ \omega} + \mbox{ a bounded term as } \omega \to 0 \,.
\end{equation}
This means that the residue of $J_{BS}(\omega)$ at $\omega=0$ is equal to 2.

To derive
\eqref{JBScoeffexp}, we write the Laurent series of $J_{BS}(\omega)$ on a punctured ball around $0$ as
$\frac{2}{\omega} + \sum_{n=0}^{\infty} a_n \omega^n$. With the notation
$\phi(z) = \sum_{n=0}^{\infty} a_n \omega^n$, we see that $\phi$ is analytic in a neighborhood of $0$ and 
\begin{equation}
\label{phidef}
J_{BS}(\omega) = \frac{2}{\omega} + \phi(\omega).
\end{equation} 
By uniqueness of analytic extension, 
$\phi$ is analytic and \eqref{phidef} holds on the domain $B_{1-\omega_1}(1)$. 
Transferring this result to the asymptotics of the coefficients $c_{J,n}$ 
yields \eqref{JBScoeffexp}.\\

(ii)
To prove the second part of the proposition, let us write $\psi(y) = J_{BS}(e^y)$. By the first part of the
proposition, $\psi$ is analytic as long as
$h$ is analytic at $e^y$. Indeed, since $e^y \neq 0$ for all complex numbers $y$, 
the only singularities of $\mathcal{J}(h (e^y))$ are due to singularities of $h$. 

The equation $e^y = \omega_1$
has infinitely many solutions of the form $i(2k+ 1) \pi + \log |\omega_1|$. 
Let us write $y_{\pm} = \pm i \pi + \log |\omega_1|$ and
$\rho_x = |y_{\pm}|$ for the singularities which are closest to the origin $y=0$. 
Then, exactly as in the proof of Proposition \ref{prop:2}, we 
see that $\psi$ is analytic on $B_{\rho_x}(0)$ and there are  
branch points
at $y_{\pm}$. The picture is now similar to that of Proposition 2.1 (ii) and is different from point (i) discussed above as the singularities closest to the origin are branch points (and are not isolated).

It remains to verify \eqref{JBScoeffexp2}. We again use the transfer results (cf. \cite{Flajolet1988}).
As in the proof of Proposition \ref{prop:2}, we define $h(\omega_1) = z_1$.
Next, we consider the Taylor expansion of $\mathcal{J}$ around $z_1 := -\eta_1^2$ (recall that $\mathcal{J}$ is analytic at $z_1$).

First, we claim that $\mathcal{J}'(z_1) =0$. Indeed, for any $\eta$ with $\eta=\tan\eta$, we have  $\eta = \frac{2\tan(\eta/2)}{1-\tan^2(\eta^2/2)}$ and so
$$
\mathcal{J}'(-\eta^2) = \frac12 - \frac{1}{2\eta} \tan\frac{\eta}{2} - \frac{1}{4\cos^2\frac{\eta}{2}}
= \frac14 - \frac{1}{2\eta} \tan\frac{\eta}{2} - \frac14 \tan^2\frac{\eta}{2} = 0$$
Thus we have the Taylor expansion
\begin{equation}\label{Jexp}
\mathcal{J}(z) = \mathcal{J}(z_1) + \frac12 \mathcal{J}''(z_1)  (z-z_1)^2
+ \frac16 \mathcal{J}'''(z_1) (z - z_1)^3 + O((z-z_1)^4)
\end{equation}
as $z \rightarrow z_1$. 
Substituting here $z - z_1 \to h(\omega) - h(\omega_1)$
and using the Puiseux expansion 
\eqref{eq:exphomega_1} gives an expansion 
of the form
\begin{equation}\label{eq:Jomegaexp}
J_{BS}(\omega) = \mathcal{J}(h(\omega)) = C_1^J (\omega-\omega_1) +
C_{3/2}^J  (\omega-\omega_1)^{3/2}
+ O((\omega - \omega_1)^2)
\end{equation}
as $\omega \rightarrow \omega_1$.
The coefficients in this expansion can be obtained in terms of the 
derivatives of $\mathcal{J}(z_1)$ and $C_i$ in \eqref{eq:exphomega_1}
by expanding $h(\omega)$ for $\omega$ close to $\omega_1$.

The branch point singularity arises from the second term in (\ref{eq:Jomegaexp})  with coefficient
\begin{equation}\label{C32J}
C_{3/2}^J = \frac16 \mathcal{J}'''(z_1) C_1^3 + \mathcal{J}''(z_1) C_1 C_2 .
\end{equation}

Finally, changing the expansion variable from $\omega-\omega_1$ to $y$ with $\omega=e^y$ is done in the same way as in the proof of Proposition \ref{prop:2}(ii).
The transfer result now says that the Taylor coefficients $A_n$ of the function $(1-z)^{3/2} = \sum_{n=0}^{\infty} A_n z^n$ satisfy
\begin{equation}
A_n =  \frac{1}{\sqrt{\pi n^5}} \left( \frac34 + O(n^{-1})\right)\, .
\end{equation}
Thus
\begin{equation}
d_{J,n} = \frac32 C_{3/2}^J \frac{(-\omega_1)^{3/2}}{\sqrt{\pi} } 
[(y_+)^{\frac32-n} + (y_-)^{\frac32-n}] n^{-\frac52} + O(n^{-\frac72}),
\end{equation}
whence \eqref{JBScoeffexp2} follows by substituting $y_\pm = \rho_x e^{\pm i\theta_x}$.

\end{proof}

The left plot in Figure \ref{Fig:J1} shows $|d_{J,n}|^{1/n}$ vs $1/n$ for the first 100 coefficients of the expansion in $y^n$, which illustrates
convergence to $1/\rho_x$, the inverse of the convergence radius in (\ref{JBSexplog}).

We study also the approximation error of the asymptotic expansion of the (signed) coefficients
\begin{equation}
\varepsilon_{J,n} = d_{J,n}/d_{J,n}^{\rm asympt} - 1 \sim O(1/n)
\end{equation}
where we denoted $d_{J,n}^{\rm asympt}$ the leading asymptotic expression for the coefficients in Eq.~\eqref{JBScoeffexp2}. 
The error is formally of order $O(n^{-1})$, so we expect that
$\varepsilon_{J,n}$ vs $1/n$ approaches zero linearly. 

The right plot in Figure \ref{Fig:J1} shows $\varepsilon_{J,n}$ vs $1/n$, which 
confirms the expected general trend towards zero of the error terms. This is seen in more
detail in the left plot in Figure \ref{Fig:J2} which zooms into the region close to the
horizontal axis. 

We note also the presence of large error outliers. 
For example for $n=97$ the relative error is $|\varepsilon_{J,97}|=5.54$
such that the subleading correction to the asymptotic
result can be comparable with the leading order coefficient.
The right plot in Fig.~\ref{Fig:J2} shows the error $\varepsilon_n$ vs 
$\cos(\theta_X(n-\frac32))$ which shows that large errors are mostly associated
with small values of the cos factor, which leads to accidental suppression in the
leading order contribution.

\begin{figure}[b!]
\centering
\includegraphics[width=3.0in]{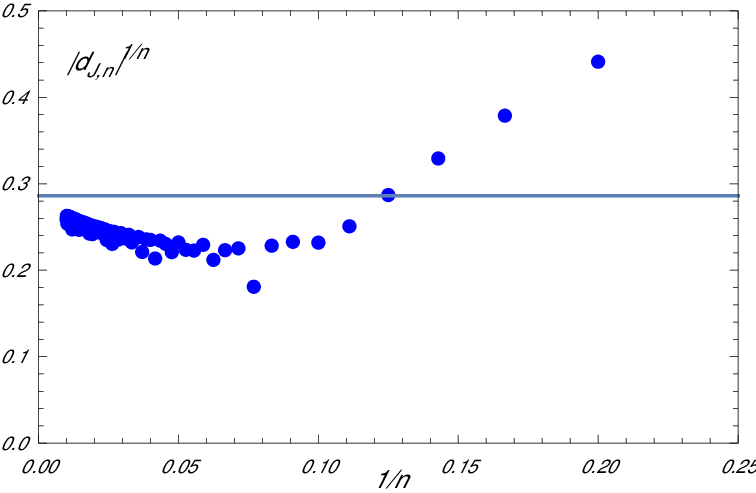}
\includegraphics[width=3.0in]{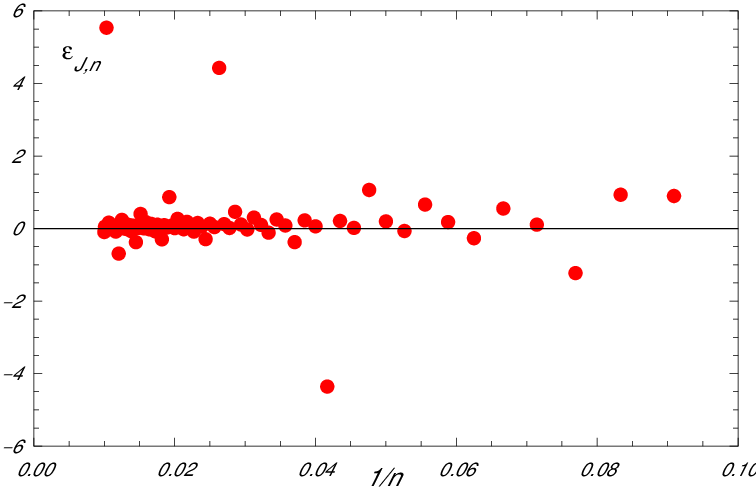}
\caption{Left: plot of $|d_{J,n}|^{1/n}$ vs $1/n$ for $n\leq 100$.
The horizontal line is at $1/\rho_x=0.286$. Right: The relative error
$\varepsilon_{J,n}=d_{J,n}/d^{\rm asympt}_{J,n}-1$ of the asymptotic coefficients
vs $1/n$.}
\label{Fig:J1}
\end{figure}

\begin{figure}[b!]
\centering
\includegraphics[width=3.0in]{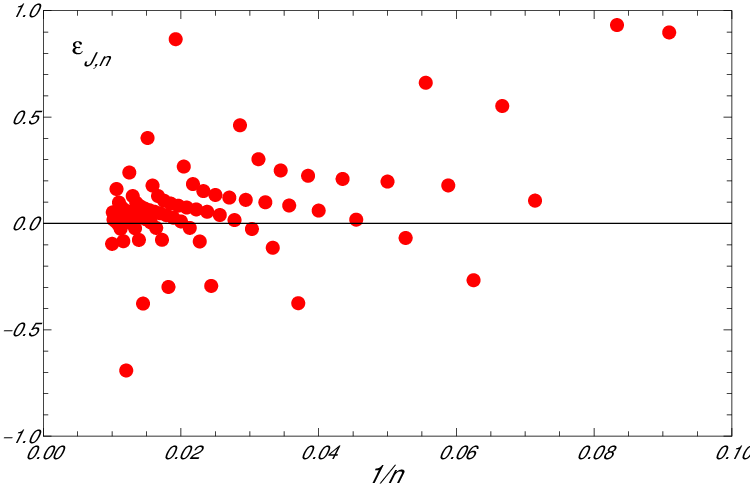}
\includegraphics[width=3.0in]{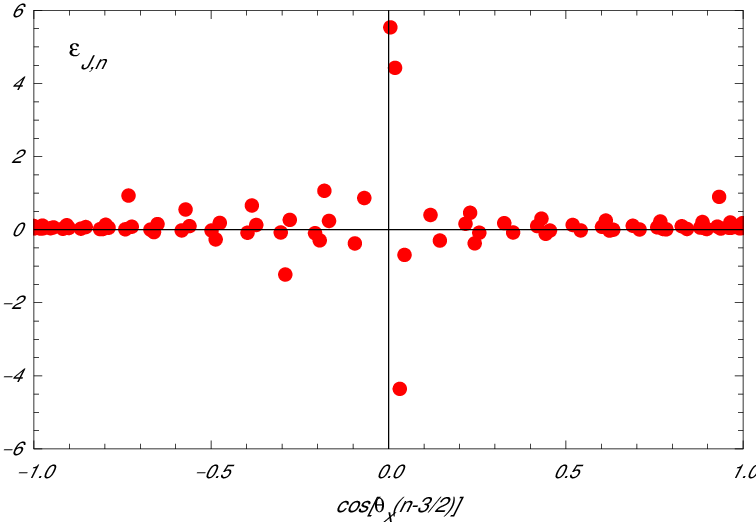}
\caption{Left: zoomed-in plot of the error of the asymptotic
coefficients $\varepsilon_{J,n}$.
Right: The relative error
$\varepsilon_{J,n}$ of the asymptotic coefficients
vs $\cos(\theta_X(n-\frac32))$.}
\label{Fig:J2}
\end{figure}


\subsection{Asymptotics of the Hartman-Watson distribution}
\label{sec:4.2}

A closed form result for the joint distribution of the time integral of the
gBM and its terminal value was obtained by Yor \cite{Yor}, see \eqref{Yor0}.

The numerical evaluation of the integral in (\ref{HWthetadef}) is challenging for 
small $t \ll 1$. Alternative ways of numerical evaluation of this distribution
have been studied in terms of analytic approximations \cite{Barrieu,Gerhold}. 

An asymptotic expansion for the Hartman-Watson distribution
$\theta_{\rho/t}(t)$ as $t\to 0$ at fixed $\rho=r  t $ was derived in \cite{HWpaper}. 
The leading order result is expressed by \eqref{HWexp0} and the real functions 
$F(\rho),G(\rho)$ are defined as
\begin{eqnarray}\label{Fsol}
F(\rho) = \left\{
\begin{array}{cc}
 \frac12 \kappa^2- \frac{\kappa}{\tanh \kappa} + \frac{\pi^2}{2}, & \, 0 < \rho < 1 \\
-\frac12 \lambda^2+ \frac{\pi - \lambda}{\tan \lambda} + \pi \lambda, & \, \rho > 1 \\
\end{array}
\right.
\end{eqnarray}
and 
\begin{eqnarray}\label{Gsol}
G(\rho) = \left\{
\begin{array}{cc}
\frac{\rho \sinh \kappa}{\sqrt{\rho \cosh \kappa - 1}}, & \,
   0 <\rho<1\\
\frac{\rho \sin \lambda}{\sqrt{1+\rho \cos \lambda }}, & \, \rho>1.\\
\end{array}
\right.
\end{eqnarray}

Here $\kappa \geq 0$ is the positive solution of the equation
\begin{equation}\label{eqx1}
\rho \frac{\sinh \kappa}{\kappa} = 1
\end{equation}
and $\lambda \in (0,\pi)$ is the solution of the equation 
\begin{equation}\label{eqy1}
\lambda + \rho\sin \lambda = \pi\,.
\end{equation}
The constraints on $\kappa$ and $\lambda$ ensure that $G(\rho)$ is positive
for positive real $\rho$.

We will study in this section the series expansions of 
$F$ and $G$ around $\rho=1$, which are relevant for the computation of the distribution (\ref{Yorlimit}) around its maximum at $a=v=1$. 

We mention that the functions $F(\rho), J_{\rm BS}(a)$ are related as 
$J_{\rm BS}(a) = 4 \inf_{\rho\geq 0} (F(\rho) + \frac{1+a^2\rho^2}{2a} - \frac{\pi^2}{2})$. This relation follows from Proposition 6 in \cite{HWpaper} and the relation $J(a) = \frac14 J_{\rm BS}(a)$. However it does not seem easy to use this relation to connect the expansions of these two functions, so we treat them separately.

The functions $\kappa(\rho),\lambda(\rho)$ are related to the function 
$h(\omega)$ studied in Section \ref{sec:2} as $\kappa(\rho)=
\sqrt{ h(1/\rho)}$
for $0 < \rho \leq 1$ and $\pi - \lambda(\rho) =i  \sqrt{ h(1/\rho)}$ for 
$\rho\geq 1$, with the usual square root function from $\mathbb R_+$ to $\mathbb R_+$.

\begin{proposition}\label{prop:FG}

(i) The function $F(e^{ - \tau})$ is expanded around $\tau=0$ as
$F(e^{-\tau}) = \sum_{n=0}^\infty d_{F,n} \tau^n$. The first few terms are
\begin{equation}\label{Fexp}
F(e^{-\tau}) = \frac{\pi^2}{2}-1 - \tau + \tau^2 + \frac{2}{15} \tau^3
+ O(\tau^4) .
\end{equation}
The large $n$ asymptotics of the coefficients is
\begin{equation}\label{dFn}
d_{F,n}  = d_F \rho_x^{-n} \cos\left( \theta_x(n-\frac32) \right) n^{-\frac52}
+ O(n^{-\frac72})
\end{equation}
with $d_F \approx -23.4047$.

ii) The function $G(e^{-\tau})$ has the expansion around $\tau=0$ as
$G(e^{-\tau}) = \sum_{n=0}^\infty d_{G,n} \tau^n$. The first few terms are
\begin{equation}\label{Gexp}
G(e^{-\tau}) = \sqrt3 \left(
1 - \frac15 \tau - \frac{1}{70} \tau^2 + \frac{1}{1050}\tau^3
+ O(\tau^4) \right).
\end{equation}
The large $n$
asymptotics of the coefficients is
\begin{equation}\label{dGn}
d_{G,n} = d_G \rho_x^{-n} \sin\left( \theta_x(n+\frac14) \right) n^{-\frac34}
+ O(n^{-\frac54})
\end{equation}
with $d_G \approx 0.719253$.

Both series expansions converge for $|\tau| < \rho_x= |\log |\omega_1| + i\pi| 
\approx 3.49295$.
\end{proposition}

\begin{proof}

(i) The function $F(\rho)$ can be written as 
{$F(\rho) = \mathcal{F}(h(1/\rho))$
with $\mathcal{F}(z)=\frac12 z - \frac{\sqrt{z}}{\tanh \sqrt{z}} + \frac{\pi^2}{2}$.} 
Now the situation (and hence the proof) is similar to Proposition~\ref{prop:J}.

The function $\mathcal{F}(z)$ has poles at {$z_k= -k^2\pi^2$}, $k\in \mathbb{Z}$
and is analytic everywhere else.
Solving the equation ${h(1/\rho)} =  -k^2\pi^2$ we find that the unique solution is $\rho = \infty$. 
Thus the function
$\tilde F (\tau) : = F(e^{- \tau})$ is analytic on $B_{\pi}(0)$ and can be extended analytically to 
$B_{\rho_x}(0)$ as in Proposition \ref{prop:J}.

In order to find the large $n$ asymptotics of $d_{F,n}$ we study again the
asymptotics of $F(e^{-y})$ around its branch points at $y_\pm$. 
This requires the Taylor
expansion of $\mathcal{F}(z)$ around $z_1=-\eta_1^2$ which is a regular point for this function.
The analysis
closely parallels that in Prop.~\ref{prop:J}(ii). 
It turns out that $\mathcal{F}'(z_1)=0$, just as for $\mathcal J(z)$
and so the expansion of
$\mathcal{F}(z)$ around $z_1$ reads
\begin{equation}
\mathcal{F}(z) = \mathcal{F}(z_1) + \frac12 \mathcal{F}''(z_1) (z-z_1)^2 +
 \frac16 \mathcal{F}'''(z_1) (z-z_1)^3  + O((z-z_1)^4).
\end{equation}
Substituting here the Puiseux expansion (\ref{eq:exphomega_1}) gives that 
the leading  singularity is a branch point at $\omega_1$ arising from the second term in 
the expansion
\begin{equation}
\mathcal{F}(h(\omega)) - \mathcal{F}(h(\omega_1)) = C_1^F (\omega - \omega_1) +
C_{3/2}^F (\omega - \omega_1)^{3/2} 
+ O((\omega-\omega_1)^2),
\end{equation}
where 
\begin{equation}
C_{3/2}^F = \frac16 \mathcal{F}'''(z_1) C_1^3 + \mathcal{F}''(z_1) C_1 C_2 \,.
\end{equation}

Changing variables to $\tau = \log\omega = -\log\rho$ gives a similar Puiseux expansion in
fractional powers of $\tau - \tau_\pm$. The asymptotics of the coefficients $d_{F,n}$ follows again by the transfer result 
\begin{equation}
d_{F,n} = \frac32 C_{3/2}^F
\frac{(-\omega_1)^{3/2}}{\sqrt{\pi} } 
[(y_+)^{\frac32-n} + (y_-)^{\frac32-n}] n^{-\frac52} + O(n^{-\frac72})
\end{equation}
and has the form (\ref{dFn}) with
\begin{equation}
d_F = \frac32
\sqrt{ \frac{(-\omega_1\rho_x)^3}{\pi } } C_{3/2}^F \approx -23.4047\,.
\end{equation}

(ii) 
The function $G(\rho)$ can be written as $G(\rho) = \mathcal{G}(h(1/\rho))$,
where we defined
\begin{equation}\label{Gcaldef}
\mathcal{G}(z):=\frac{\sqrt{z}}{\sqrt{\frac{\sqrt{z}}{\tanh \sqrt{z}}-1}} \,.
\end{equation}
Writing the numerator in this form ensures that $G(\rho)$ is positive for real and 
positive $\rho$, which is the case required for the application to the asymptotics 
of the Hartman-Watson distribution (\ref{HWexp0}).
(This condition was ensured in the original expression by imposing the constraints $\kappa \geq 0$ and $\lambda\in (0,\pi)$.) 
This condition replaces the even property of
$\mathcal{F}(z)$ which was used to resolve the sign ambiguity before.

The function $\mathcal{G}(z)$ has branch points on the negative real axis at
$z_k= - \eta_k^2$ where the denominator in (\ref{Gcaldef}) vanishes.
The branch point closest to the origin $z=0$ is at $z_1$ and will be relevant
for our purpose
as $z_1$ is mapped to $\omega_1 = 1/\rho$, the branch point of
$h$.

In the neighborhood of $z_1$ the function $\mathcal{G}(z)$ diverges as
\begin{equation}\label{calGexp}
\mathcal{G}(z) = \sqrt{2z_1} (z-z_1)^{-\frac12} + {O((z-z_1)^{\frac12})} \,.
\end{equation}
The asymptotics of $\mathcal{G}(h(1/\rho))$ around $1/\rho=\omega_1$ are
obtained by substituting the Puiseux series (\ref{eq:exphomega_1}) into (\ref{calGexp}). Denoting
$x=1/\rho$ for simplicity, we have as $x\to \omega_1$
\begin{equation}
\mathcal{G}(h(x)) = \sqrt{\frac{2z_1}{ C_1}}  (x -\omega_1)^{-\frac14} + {O((x-\omega_1)^{\frac14})}\,.
\end{equation}

This asymptotics can be transferred to the large $n$ asymptotics of the
coefficients $d_{G,n}$ as
\begin{equation}
d_{G,n}  = \frac{1}{\Gamma(1/4)n^{3/4}}
\sqrt{\frac{-2\eta_1^2}{C_1}} 
(-\omega_1)^{-\frac14} 
[(y_+)^{-\frac14-n} - (y_-)^{-\frac14-n}] + O(n^{-5/4})\,.
\end{equation}
We chose the branches of the $(y-y_\pm)^{-1/4}$ functions such that
they add up to a real function along the positive real axis. 
This introduces a minus sign between the two terms.
This reproduces the stated result, with 
\begin{equation}
d_G = \frac{1}{\Gamma(1/4)}
\sqrt{\frac{2\eta_1^2}{C_1}} 
(-\omega_1\rho_x)^{-\frac14} \approx 0.719253 \,.
\end{equation}

\end{proof}

We present next some numerical tests of the asymptotic expansions for $F(\rho)$ and $G(\rho)$. Figures \ref{Fig:F1} and \ref{Fig:G1} illustrate the performance of the
asymptotic expansions for these two functions. They are similar, so we comment in detail only the plots for $F(\rho)$.

The left plot in Figure \ref{Fig:F1}  shows $|d_{F,n}|^{1/n}$ vs $1/n$, which shows 
good convergence to $1/\rho_x$, which is the inverse of the convergence radius for the expansion (\ref{Fexp}). The right plot shows the approximation error of the asymptotic expansion of the coefficients
$\varepsilon_{F,n} = d_{F,n}/d_{F,n}^{\rm asympt} -1 \sim O(1/n)$ vs $1/n$,
where 
$d_{F,n}^{\rm asympt}$ is the first term on the right hand side of
\eqref{dFn}.
The numerical results confirm a general decreasing trend of the approximation error 
with $n$ as expected. 

\begin{figure}[b!]
\centering
\includegraphics[width=3.0in]{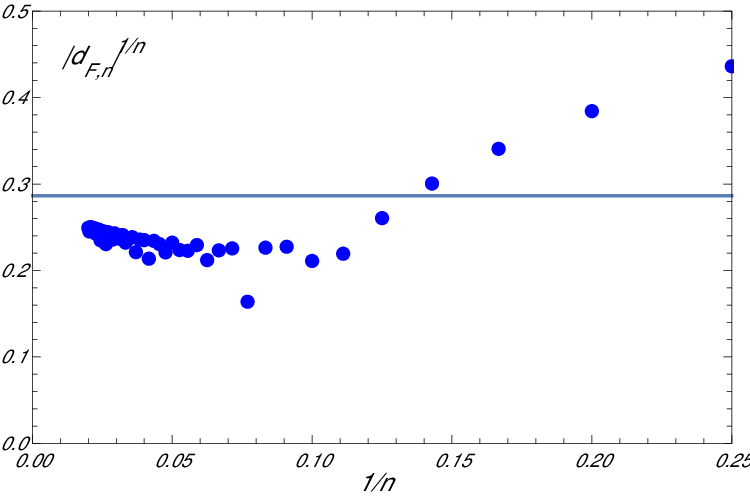}
\includegraphics[width=3.0in]{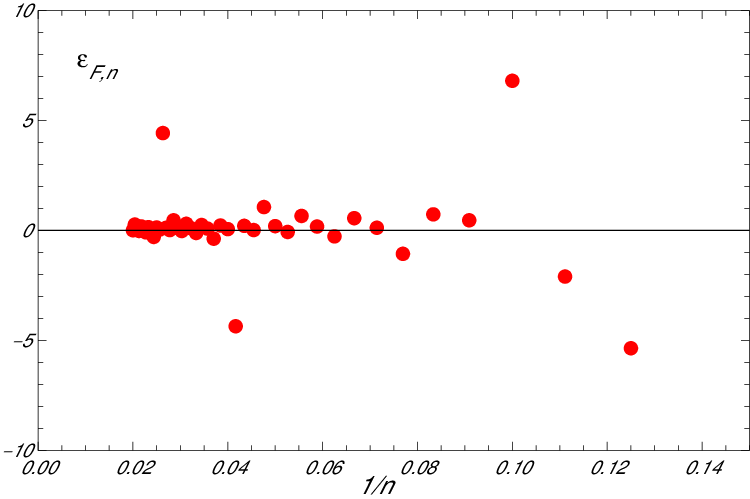}
\caption{Left: plot of $|d_{F,n}|^{1/n}$ vs $1/n$ for $n\leq 100$.
The horizontal line is at $1/\rho_x=0.286$. Right: The relative error
$d_{F,n}/d^{\rm asympt}_{F,n} -1$ of the asymptotic coefficients
vs $1/n$.}
\label{Fig:F1}
\end{figure}

\begin{figure}[b!]
\centering
\includegraphics[width=3.0in]{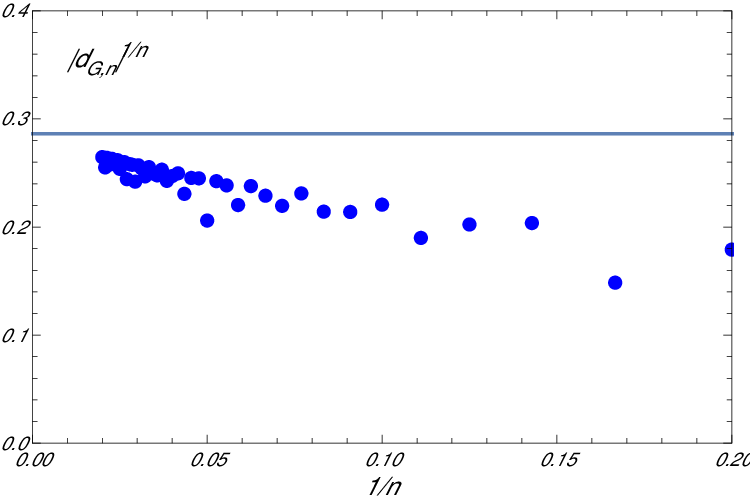}
\includegraphics[width=3.0in]{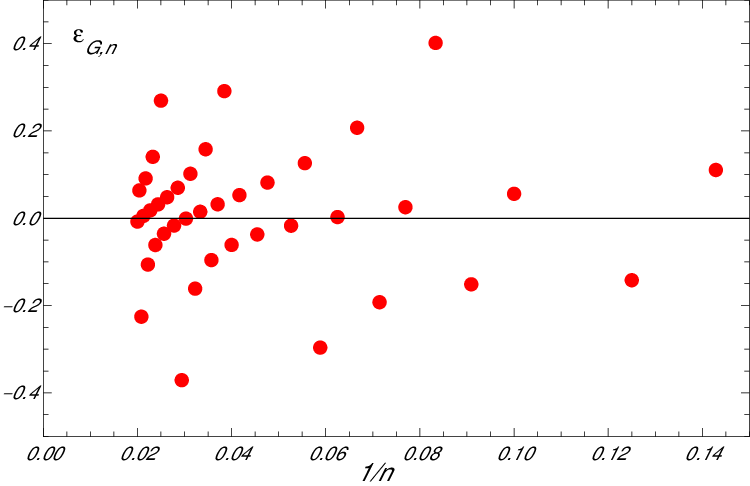}
\caption{Left: plot of $|d_{G,n}|^{1/n}$ vs $1/n$ (blue points) and the
asymptotic coefficients $|d^{\rm asympt}_{G,n}|^{1/n}$ (red points) vs $1/n$
for $n\leq 50$.
The horizontal line is at $1/\rho_x=0.286$. Right: The relative error
$d_{G,n}/d^{\rm asympt}_{G,n} -1$ of the asymptotic coefficients
vs $1/n$.}
\label{Fig:G1}
\end{figure}

\clearpage
\section{Application: Asian options pricing in the Black-Scholes model}
\label{sec:5}

Here, we provide a simple application of the expansions studied in this
paper, namely the pricing of an Asian option in the Black-Scholes model.
Assume that the asset price follows a geometric Brownian motion
\begin{equation}
dS_t = r S_t dt  + \sigma S_t dW_t
\end{equation}
started at $S_0$,
or equivalently, by It\^o's formula, 
$S_t = S_0 e^{(r - \sigma^2 /2)t + \sigma W_t}$.
Of course, in practice one has to consider the issue of 
mis-specification of the 
price process, which introduces model risk.  

Asian option prices are given by expectations in the risk-neutral measure
\begin{eqnarray}
&& C_A(K,T) = e^{-rT} \mathbb{E}
\Big[\Big(\frac{1}{T} \int_0^T S_t dt - K\Big)^+\Big]\,,\\
&& P_A(K,T) = e^{-rT} \mathbb{E}\Big[\Big(K - \frac{1}{T} \int_0^T S_t dt\Big)^+\Big]\,.
\end{eqnarray}
$C$ and $P$ stand for call and put options and the subscript $A$ denotes Asian.

Using the self-similar property of the Brownian motion, the Asian option prices can be expressed in terms of expectations of the standardized integral of gBM 
$A_t^{(\mu)}=\int_0^t e^{2(B_s+2\mu s)} ds$ as 
\begin{equation}
C_A(K,T) = e^{-rT} S_0 c_A(k,\tau)\,,\quad
P_A(K,T) = e^{-rT} S_0 p_A(k,\tau)
\end{equation}
where
\begin{equation}\label{cpA}
c_A(k,\tau) := \mathbb{E}\Big[\Big(\frac{1}{\tau} A_\tau^{(\mu)} - k\Big)^+\Big]\,,\quad
p_A(k,\tau) := \mathbb{E}\Big[\Big(k - \frac{1}{\tau} A_\tau^{(\mu)} \Big)^+\Big]\,.
\end{equation}
and the expectations in (\ref{cpA}) are expressed in terms of reduced parameters
\begin{equation}
\tau = \frac14 \sigma^2 T\,,\qquad
\mu =\frac{2r}{\sigma^2} - 1 \,,\qquad
k = \frac{K}{S_0}\,,
\end{equation}
see e.g. \cite{Geman1993}.
Denote the density of the standardized time-average of the gBM
\begin{equation}
\mathbb{P}\Big(\frac{1}{t} A_t^{(\mu)} \in da \Big) = f(a,t) \frac{da}{a} \,.
\end{equation}
The leading $t\to 0$ asymptotics (\ref{Yorlimit}) gives an approximation for this 
density
\begin{equation}
f(a,t) = f_0(a,t) (1 + O(t))\,,\quad
f_0(a,t) := \frac{1}{n(t)} \frac{1}{2\pi t} e^{-\frac12\mu^2 t}
a^\mu \int_0^\infty \rho^\mu G(\rho) e^{-\frac{1}{t} I(a,a \rho)} \frac{d\rho}{\rho}
\end{equation}
where $n(t)$ is a normalization integral, defined such as to ensure the normalization condition $\int_0^\infty f_0(a,t) \frac{da}{a} = 1$ for all $t>0$.

Recalling \eqref{eq:FI}, the expectations (\ref{cpA}) can be approximated as integrals 
over the $f_0(a,t)$ density as
\begin{eqnarray}\label{cAint}
&& c_A^{(0)}(k,\tau) = \int_k^\infty (a-k) f_0(a,t) da = \\
&& =\frac{1}{n(\tau)} \frac{1}{2\pi t} e^{-\frac12 \mu^2 t} \int_0^\infty \rho^\mu G(\rho) e^{-\frac{1}{t}[F(\rho)-\frac{\pi^2}{2}]}
\left( \int_k^\infty a^\mu (a-k) e^{-\frac{1+a^2 \rho^2}{2a t}} \frac{da}{a} \right)
\frac{d\rho}{\rho} \nonumber
\end{eqnarray}
and analogously for $p_A(k,\tau)$.

The normalization factor $n(\tau)$ can be expressed as a one-dimensional integral
\begin{equation}\label{nint}
n(\tau) = \frac{1}{\pi \tau} e^{-\frac12\mu^2\tau} \int_0^\infty G(\rho) K_{-\mu}(\rho/\tau) e^{-\frac{1}{\tau} (F(\rho) -\frac{\pi^2}{2})} \frac{d\rho}{\rho}\,.
\end{equation}

\subsection{Details of implementation}

For the purpose of evaluation of the integrals in (\ref{cAint}) we 
represent the functions $F(\rho),G(\rho)$ as follows: 

i) within the convergence
domain of their expansions (\ref{Fexp}) and (\ref{Gexp}) we use the 
series expansions derived in Proposition \ref{prop:FG},  truncated to an appropriate order $N_F, N_G$.

ii) outside the convergence domain we use the asymptotic expansions derived in points (i),(ii) Propositions 4 and 5 of \cite{HWpaper} for $F(\rho),G(\rho)$ as $\rho\to 0$ and $\rho\to \infty$, respectively.

This gives the following approximation for the functions $F(\rho),G(\rho)$, defined by
\begin{equation}
\bar F(\rho) := \left\{
\begin{array}{cc}
F_L(\rho) & \log\rho < \rho_L \\
\sum_{k=0}^{N_F} d_{F,k} \log^k(1/\rho) & \log \rho \in [\rho_L,\rho_R] \\
F_R(\rho) & \log \rho > \rho_R \\
\end{array}
\right.
\end{equation}
and
\begin{equation}
\bar G(\rho) := \left\{
\begin{array}{cc}
G_L(\rho) & \log\rho < \rho_L \\
\sqrt3 \sum_{k=0}^{N_G} d_{G,k} \log^k(1/\rho) & \log \rho \in [\rho_L,\rho_R] \\
G_R(\rho) & \log \rho > \rho_R \\
\end{array}
\right.
\end{equation}
where $[\rho_L,\rho_R] \subset [-\rho_x,\rho_x]$ with $\rho_x = 3.49295$ the
convergence radius of the series expansions. $F_L(\rho), G_L(\rho)$ denote the
$\rho \to 0$ asymptotic expansions given in equations (49) and (58) of \cite{HWpaper} and $F_R(\rho), G_R(\rho)$ denote the $\rho\to \infty$ asymptotic expansions given in equations (50), (59) of \cite{HWpaper}. 
The first few coefficients $d_{F,k},d_{G,k}$ are tabulated in Table~\ref{tab:dFG}.
The approximation $\bar F(\rho), \bar G(\rho)$ can be made arbitrarily precise by increasing the truncation orders $N_{F,G}$ and by including higher order terms in the tail asymptotics.

\begin{table}[htbp]
\centering
\caption{The first 10 coefficients of the series expansions of $F(\rho), G(\rho)$ in powers of $\log(1/\rho)$. 
The leading terms of these series expansions are $d_{F,0}=\frac{\pi^2}{2}-1$
and $d_{G,0}=1$.}
\begin{tabular}{|c|c|c||c|c|c|}
\hline
$k$ & $d_{F,k}$ & $d_{G,k}$ & $k$ & $d_{F,k}$ & $d_{G,k}$\\
\hline
\hline
1 & 1 & $ \frac15$ & 6 & $\frac{4742}{3,031,875}$ 
   & $-\frac{107,749}{10,032,750,000}$ \\
2 & 1 & $-\frac{1}{70}$ & 7 & $-\frac{43,636}{197,071,875}$ 
   & $\frac{27,333,619}{1,876,124,250,000}$ \\
3 & $-\frac{2}{15}$ & $ - \frac{1}{1050}$ & 8 & $\frac{146,287}{6,897,515,625}$ 
   & $-\frac{308,907,281,743}{109,790,791,110,000,000}$ \\
4 & $\frac{19}{525}$ & $\frac{299}{323,400}$ & 9 
      & $-\frac{68,146}{57,984,609,375}$ 
      & $-\frac{1,589,498,602,063}{4,940,585,599,950,000,000}$ \\
5 & $-\frac{22}{2,625}$ & $-\frac{96,917}{525,525,000} $
   & 10 & $\frac{6,740,719,066}{38,598,324,999,609,375}$ 
   & $\frac{28,340,195,926,465,733}{103,406,456,606,953,500,000,000}$ \\
    \hline
    \end{tabular}%
  \label{tab:dFG}%
\end{table}%

\subsection{Numerical examples}

We illustrate the application of the method on the seven test cases introduced in \cite{Fu1998} and commonly used in the literature as benchmark tests for Asian options pricing. 
Table \ref{tab:7scenarios} shows the option prices following from the method 
proposed here for each of these scenarios, comparing with the precise evaluation 
in Linetsky \cite{Linetsky} obtained using the spectral expansion method.

\begin{table}[htbp]
\centering
\caption{Numerical tests for pricing Asian options in the Black-Scholes model on
the seven scenarios commonly considered in the literature \cite{Fu1998,Linetsky}.
All scenarios have $K=2.0$. The results from the method proposed here are shown in the column $C_A(K,T)$ which are compared with the precise benchmarks from the spectral expansion method from Linetsky \cite{Linetsky}.}
\begin{tabular}{c|cccc|cc|c||c|c|c}
\hline
Scenario & $S_0$ & $r$ & $\sigma$ & $T$ 
               & $\mu$ & $\tau$ & spectral \cite{Linetsky} 
               & $c_A(k,\tau)$ & $n(\tau)$ & $C_A(K,T)$ \\
\hline
\hline
1 & 2.0 & 0.02 & 0.10 & 1 
   & 3 & 0.0025 & 0.055986 
   & 0.028543 & 1.00004 & 0.055954 \\
2 & 2.0 & 0.18 & 0.30 & 1 
   & 3 & 0.0225 & 0.218387 
   & 0.130771 & 1.00032 & 0.218388 \\
3 & 2.0 & 0.0125 & 0.25 & 2 
   & -0.6 & 0.03125 & 0.172269 
   & 0.088354 & 1.00045 & 0.172269 \\
4 & 1.9 & 0.05 & 0.50 & 1
   & -0.6 & 0.0625 & 0.193174 
   &  0.106978 & 1.00089 & 0.193174 \\
5 & 2.0 & 0.05 & 0.50 & 1 
   & -0.6 & 0.0625 & 0.246416 
   & 0.12964 & 1.00089 & 0.246415 \\
6 & 2.1 & 0.05 & 0.50 & 1 
   & -0.6 & 0.0625 & 0.306220 
   & 0.153432 & 1.00089 & 0.306220 \\
7 & 2.0 & 0.05 & 0.50 & 2 
   & -0.6 & 0.125 & 0.350095 
   & 0.193799 & 1.00177 & 0.350093  \\
    \hline
    \end{tabular}%
  \label{tab:7scenarios}%
\end{table}%

We used an approximation for $\bar F(\rho), \bar G(\rho)$ keeping $N_F=N_G=6$
terms in the series expansion. Adding more terms has no impact on the numerical
values shown. The series expansion was used in the range $\rho\in [0.04,32.88]$ which is included in the convergence domain of the series expansions. The impact of the tails region on the results in Table \ref{tab:7scenarios}  is negligible.

For the numerical evaluation we changed the integration variable $\rho$ to 
$z = \log\rho$ in both integrals in (\ref{cAint}) and (\ref{nint}).  
The 2-dimensional integration was performed in \textit{Mathematica}
using \texttt{NIntegrate} with \texttt{Method -> NewtonCotesRule} and \texttt{MaxIterations -> 100}.

The results in Table \ref{tab:7scenarios} show good agreement with the precise benchmark values of \cite{Linetsky}, the difference being less than 1\% in all cases. 
For most scenarios, the impact of including terms beyond quadratic order in the joint distribution (\ref{Yorlimit}) is larger than the impact of varying the volatility parameter by 
$\Delta \sigma \sim 0.2\%$, which makes them relevant in practical applications where $\sigma$ is known with precision of this order of magnitude.
\vspace{0.5cm}

\textbf{Acknowledgments.}
We thank two anonymous referees for helpful comments and suggestions. 
The research of PN was partially sponsored by NSF DMS 1952876 and the Charles Simonyi Endowment at the Institute for Advanced Study, Princeton, NJ.

\end{document}